\documentstyle{amsart}
\def\Om{{ \Omega }}
\def\ra{{ \rightarrow }}
\def\om{{ \omega }}
\newtheorem{thm}{Theorem}

\theoremstyle{definition}

\begin{document}

\title[Diophantine approximation and deformation]
{Diophantine approximation and deformation}
\author{Minhyong Kim, Dinesh S. Thakur and Jos\'e Felipe Voloch}
\address{University of Arizona, Tucson and University of Texas, Austin}
\thanks{Authors supported in part by NSF and NSA grants}
\maketitle

\vskip .2truein

\abstractname{: It is well-known that while the analogue of
Liouville's theorem on diophantine approximation holds in finite
characteristic, the analogue of Roth's theorem fails quite badly. We
associate certain curves over function fields to given algebraic power
series and show that bounds on the rank of Kodaira-Spencer map of this
curves imply bounds on the diophantine approximation 
exponents of the power series, with more
`generic' curves (in the deformation sense) giving lower exponents.
If we transport Vojta's conjecture on height inequality to finite
characteristic by modifying it by adding suitable deformation
theoretic condition, then we see that the numbers giving rise to
general curves approach Roth's bound. We also prove a hierarchy of 
exponent bounds for approximation by algebraic quantities of 
bounded degree. }

\vskip .2truein

\setcounter{section}{-1}
\section{Diophantine approximation exponents}

0.0 Let $F$ be a finite field of characteristic $p$. 
For $\beta$ an element of $F((t^{-1}))$ algebraic irrational over $F(t)$
(an algebraic irrational real number, respectively), define
its diophantine approximation exponent $E(\beta)$ by
$$E(\beta):=\lim\sup (-\frac{\log |\beta-P/Q|}{\log |Q|})$$
where $P$ and $Q$ run over polynomials in 
$F[t]$ (integers, respectively), the absolute value is the usual 
one in each case and the limit is taken as $|Q|$ 
grows. 

0.1 The well-known theorems of Dirichlet and Liouville 
and their analogues for function fields [M] show 
that $2\leq E(\beta)\leq d(\beta)$, where $d(\beta)$ 
is the algebraic degree of $\beta$ defined as $[F(t,\beta):F(t)]$
($[{\mathbb Q}(\beta):{\mathbb Q}]$, respectively). That the 
diophantine approximation results and in particular, the improvement
on the Liouville bound of $d(\beta)$ have interesting implications 
for the study of related diophantine equations is well-known since the 
work of Thue, Siegel etc. For the real 
number case, the well-known theorem of Roth shows 
that $E(\beta)=2$, but Mahler showed [M] that 
$E(\beta)=d(\beta)=q$ for $\beta=\sum t^{-q^i}$, 
(so that $\beta^q-\beta-t^{-1}=0$)
as  a straightforward estimate of approximation by truncation 
of this series shows. (Here and in what follows $q$ is a 
power of $p$). 
Osgood [O2] and Baum and Sweet [BS] gave many examples 
in various degrees. See [T] for the references to other examples. 

0.2 For given $d=d(\beta)$, $E(\beta)$'s form a countable subset of
interval $2\leq x\leq d$. What is it?  Does it contain any irrational
number? Does it contain all the rationals in the range?  
In [Sc2] and [T], the following
result was proved.

\begin{thm} 
 Given any rational $\mu$ between $2$ 
and $q+1$,  we can find a family of $\beta$'s 
(given by explicit equations and explicit 
continued fractions),  with $E(\beta)=\mu$ and $d(\beta)
\leq q+1$. 
\end{thm}

The question of exact degree of
$\beta$  is easily addressed for explicit families.
(For more, see [T]). 

\section{Differential equations and deformations}

1.0 Osgood [O2] proved  that the Liouville/Mahler bound can be
improved to (even effectively) $E(\beta)\leq 
\lfloor (d(\beta)+3)/2\rfloor$ 
(or rather $\lfloor d(\beta)/2\rfloor +1$, see [Sc1] or [LdM1])  
for
$\beta$ not satisfying the generalized Riccati differential equation
$d\beta/dt=a\beta^2+b\beta +c$, with $a, b, c\in F(t)$.
Most known (see e.g., [L] for exceptions) examples $\beta$, whose
continued fraction is known, do satisfy Riccati equation and indeed
$\beta$ is an integral linear fractional transformation 
of $\beta^q$.
What is the range of exponents for $\beta$ not of this
form? 

1.1 In [V2],  [V3] there is an observation 
(for the lack of a better reference, we provide a proof of 
this in the appendix) that 
the Riccati condition is equivalent to the vanishing of 
Kodaira-Spencer (we write KS in  short-form) 
class of projective line minus conjugates 
of $\beta$. Hence, it may not be too wild to speculate 
that it might be possible to successively improve on Osgood's bound, if
we throw out some further classes of differential 
equations coming from the conditions that some  corresponding
Kodaira-Spencer map (or say the vector space 
generated by derivatives of the cross-ratios of conjugates of 
$\beta$) has rank not more than some integer. 
It should be also noted that even though the KS
connection holds in characteristic zero, analogue 
of Roth's theorem  holds in the complex function field case. 
The Osgood bound 
still holds (conjectured in [V1], proved in [LdM1] and again in 
[LdM2]) 
by throwing out only a subclass given by `Frobenius' 
equation $\beta^q=(a\beta+b)/(c\beta +d)$. 
This might be the best one can get.
Similarly, the differential equation hierarchy suggested above 
might have some corresponding more refined hierarchy.

\section{Height inequalities for algebraic points}

2.0 Though we have not succeeded yet in improving the Osgood 
bound unconditionally, by throwing out more classes of numbers
(we want good conditions like Osgood rather than a trivial 
way of obtaining this by throwing solutions of 
$y'=P_k(y)$, where $P_k$ is a polynomial of degree
$k$, for $j\leq k\leq d$ to force  exponents to be 
less than $j$ via the Kolchin theorem [O1] on denomination, 
combined with the proposition of pa. 762 of [Sc1]),
we prove  existence of  hierarchies, given by deformation 
theoretic conditions,  of 
bounds, using results of [K], which we now recall: 
Let $X$ be a smooth 
projective surface over a perfect field $k$. Assume that 
$X$ admits a map $f: X\rightarrow S$ to a smooth projective 
curve $S$ defined over $k$, with function field $L$ in such 
a way that the fibers of $f$ are geometrically connected curves and the 
generic fiber $X_L$ is smooth of genus $g\geq 2$. Consider 
algebraic points $P: T\rightarrow X$ of $X_L$, where 
$T$ is a smooth projective curve mapping to $S$ (such that the 
triangle commutes). 
Define the canonical height of $P$ to be $h(P)
:=\deg P^*\omega/[T:S]=\langle P(T).\omega\rangle/[K(P(T)):L]$, where 
$\omega=\omega_X:=K_X\otimes f^*K_S^{-1}$ denotes the 
relative dualising sheaf for $X\rightarrow S$. This is a 
representative for the class of height functions on $X_L(\overline{L})$
associated to the canonical sheaf $K_{X_L}$. Define the relative 
discriminant  to be $d(P):=(2g(T)-2)/[T:S]=(2g(P(T))-2)/[K(P(T)):L]$.

The Kodaira-Spencer map is constructed on any open set $U\subset S$
over which $f$ is smooth from the exact sequence 
$0\rightarrow f^*\Omega_U^1 \rightarrow \Omega_{X_U}^1\rightarrow
\Omega_{X_U/U}^1\rightarrow 0$, by taking the coboundary map 
$KS: f_*(\Omega_{X_U/U}^1)\rightarrow \Omega_U^1\otimes R^1f_*(
{\cal O}_{X_U})$.

\begin{thm} 
{\mbox [K]} (1) Suppose the KS map of $X/S$
(defined on some open subset of $S$) is non-zero. 
Then 
$h(P)\leq (2g-2)d(P)+O(h(P)^{1/2})$ if 
$g>2$. If g=2, then $h(P)\leq (2+\epsilon)d(P)+O(1)$.

(2) Suppose the KS map of $X/S$ 
has maximal rank, then $h(P)\leq (2+\epsilon)d(P)+O(1)$.
\end{thm}

2.1  The inequality in (2) was proved [Voj2] in the
 characteristic $0$ function field analogue, without any hypothesis,
 by Vojta, who also conjectured [Voj1] the stronger inequality with $2$
 replaced by $1$ in the number field (and presumably also in the
 characteristic $0$ function field) case.

2.2  Modifying  the proof in [K] of Theorem 2, we get 
an hierarchy of bounds: 

{\bf Claim:} If the rank of the  kernel of the KS map is $\leq i$, then 
we have $$h(P)\leq (Max((2g-2)/(g-i), 2)+\epsilon)d(P) + O(1),\ \ (0\leq
i<g)$$ (Note that the maximum is $2$ only for $i=0, 1$. )

{\bf Proof:} To be consistent with the notation of 
[K], we change notation in 2.2 only: let $F$ to be the function field of 
$S$ and  $L$ to be  a line bundle. The claim  follows 
by combining the last displayed inequality in the proof of 
Theorem 1 of [K] with the argument in the proof 
of Theorem 2 connecting $\deg (G_F)$ to the rank of the 
kernel of the KS map. 

In more detail, in [K], we have constructed a finite collection of exact
sequences $0\ra L \ra \Omega_X\ra G\ra 0$ such that all points $P:T\ra
X$ not satisfying $h(P)\leq (2+\epsilon)d(P)+O(1)$ will be tangential
to some $L$, i.e., the composed map $0\ra L \ra \Om_X \ra P_*(\Om_T)$
will be zero, which implies
that there is a non-zero map $P^*G \ra \Om_T$ giving us an inequality
$h_G(P) \leq d(P)$ for the height with respect to $G$. The bound
$[(2g-2)/(g-i)+\epsilon]d(P)+O(1)$ for the canonical height follows
from comparing the two heights using a lower bound for the degree of
$G_F$.  That is,$$\deg G_F=2g-2-h^0(L_F)+h^1(L_F)-g+1 \geq
g-1-i+h^1(L_F)\geq g-i.$$ This follows from two facts: First, consider
the exact sequence
$$0\ra \Om_F\ra H^0((\Om_X)_F)\ra 
H^0((\om_X)_F)\ra H^1({\cal O}_{X_F})\otimes \Om_F$$
appearing in the definition of the KS map (the last arrow).
Any subspace of $H^0((\Om_X)_F)$ not contained in $\Om_F$
contributes to the kernel of KS. Now, if $\deg L_F\leq 0$,
then $\deg G_F\geq 2g-2$. So we may assume $\deg L_F >0$.
But, then $L_F$ is not contained in $f^*\Om_F$ (which has degree zero)
and hence intersects with it trivially (because both are saturated
subsheaves of $(\Om_X)_F$). Thus, $H^0(L_F)\cap \Om_F=0$
and $H^0(L_F)$ injects into the kernel of the KS map.
Thus, we get $H^0(L_F)\leq i$ by assumption. On the other
hand, since $L_F$ is not contained in $f^*\Om_F$,
it must possess a non-trivial map to $\om_{X_F}$
from which we get that $L_F$ is special, i.e.,
$h^1(L_F)\geq 1$.

For heights, we therefore get
$$h(P) \leq [(2g-2)/(g-i) +\epsilon] 
h_G(P)+O(1)\leq [(2g-2)/(g-i) +\epsilon]d(P)+O(1)$$

 More precise knowledge of jumps in 
the bounds would  depend on the fine gap structure for $G_F$.

2.3 It would be of interest to have a nice geometric condition that would
allow us to extend to all points the bounds we get for degenerate
points, since this would give us better than a $2+\epsilon$ bound
(namely, $2-2/g+\epsilon$) in the maximal rank case. But the argument
for general points at the beginning of [K] cannot be improved with the
techniques of that paper. We also do not know whether the points 
we use below are degenerate or not. 

\section{Exponent hierarchy}

3.0 We now apply this theorem to get bounds 
on the exponents for  the diophantine approximation 
situation by associating to $\beta=\beta(t)$ some curves $X$ 
over $F(t)$ and 
associating to its approximations some algebraic points 
$P$ on them.

Let $f(x)=\sum_{i=0}^d f_ix^i$ be an irreducible polynomial with
$\beta=\beta(t)$ as a root and with $f_i\in F[t]$ being relatively
prime. Let $F(x, y)= y^df(x/y)$ be its homogenization. (there 
should be no confusion with the field $F$). 

3.1  Assume $p$ does not divide $d$. 
Let $X$ be the (projective) Thue curve with its affine 
equation $F(x, y)=1$.  Given a rational approximation 
$x/y$ (reduced in the sense that $x, y\in F[t]$ are relatively 
prime) to $\beta$, with $F(x, y)=m(t)$, we associate the
algebraic point $P=(x/m^{1/d}, y/m^{1/d})$ of $X$.
Both $F_x$ and $F_y$ are not simultaneously zero 
on $F=1$ by Euler's theorem on homogeneous functions and at 
infinity there are $d$ distinct points given by 
$F(x, 1)=0$, so that $X$ is non-singular and Theorem 2 
can be applied. Then $g=(d-1)(d-2)/2$, so that 
$2g-2=d^2-3d$, and $[K(P(T)):L]=d$. Note that $g\geq 2$ implies 
that $d\geq 4$, which implies $g\geq 3$.

Now the naive height of $P$ is $\deg(y)$
(if $\deg(x)$ is bigger than $\deg(y)$, it differs by the fixed
$\deg(\beta)$, so the difference does not matter below, similarly we
ignore $\epsilon$'s  which do not matter at the end). Hence
$h(P)+O(h(P)^{1/2})=(2g-2)\deg(y)/d=(d-3)\deg(y)$.

Let $e=E(\beta)$. We want upper bounds for $e$. If $x/y$ is an
approximation approaching the exponent bound, then degree of the
polynomial $m(t)$ is asymptotically (as $\deg(y)$ tends to infinity)
$(d-e)\deg (y)$.  Since $K(P)$ over $L$ is totally ramified at zeros
of $m$ and at infinity, by the Riemann-Hurwitz formula, we have (since
$p$ does not divide $d$) $2g(P)-2\leq -2d +((d-e)\deg(y)+1)(d-1)$,
which is asymptotic to $(d-e)(d-1)\deg(y)$.

Hence under the maximal rank hypothesis, the Theorem gives us $d-3\leq
2(d-1)(d-e)/d$. This simplifies to $e\leq d/2+ d/(d-1)$, which is
slightly worse than Osgood bound, but approaches it for large even $d$. 

3.1.2 Also note that if Vojta's conjectured inequality is assumed to hold in
characteristic $p$ under the maximal rank hypothesis (this may be
reasonable to do, taking into account parallel results (Remark 1
above) in the two cases), we get $e\leq 2d/(d-1)$, which tends to the
Roth bound $2$ as $d$ tends to infinity.

3.2  Let $X$ have affine equation  $y^k=f(x)$, with 
 $k$ relatively prime with $p$ and $d$.  
Corresponding to a (reduced) approximation 
$x/z$, let $P=(x/z, (m/z^d)^{1/k})$, where $m=F(x, z)$. 
Then by Riemann-Hurwitz (as $p$ does not divide $k$), we have 
$g=(d-1)(k-1)/2$, so that $2g-2=(d-1)(k-1)-2$,  and $[K(P(T)):L]=k$. 
We assume that $g>1$. 

Now  the naive height is $\deg (z)$ ($\deg(x)$ differs by
an additive constant, so it does not matter which is bigger). Hence
the $h(P)+ O(h(P)^{1/2})$ is $((d-1)(k-1)-2))/k$ times that. 

Now zeros of $m$ and $z$ can ramify totally, so that Riemann-Hurwitz
(as $p$ does not divide $k$) gives (for approximation approaching the
exponent bound) $2g(P)-2\leq -2k+(1+d-e)(k-1)\deg (z)$, which is
$(1+d-e)(k-1)\deg(z)$ asymptotically.  Hence under the maximal rank
hypothesis, the Theorem gives $(d-1)(k-1)-2\leq 2(1+d-e)(k-1)$, which
simplifies to $e\leq (d+3)/2 +1/(k-1)$. This is again worse than, but
asymptotic to, Osgood bound.

3.2.2 In this case, if we assume Vojta's bound under the maximal rank
hypothesis, then we get $e\leq 2+2/(k-1)$ again approaching Roth
bound, this time with $k$ approaching infinity.  So we can say that
$e=2$ under the maximal rank hypothesis and assuming the corresponding
modification of Vojta's conjecture.

\section{Approximation by algebraic functions of bounded degree}

4.0 Similar ideas can be used to study approximation of $\beta$ by
algebraic functions of lower degree in the spirit of Wirsing's theorem
[Sc3]. The setting is as follows: Let $\beta$ as before of degree $d$
over $F(t)$. Now we want to see how close $\beta$ can be to $\alpha$
of degree $r < d$ over $F(t)$. Let $\alpha$ have height $H$ and be
such that $-\log|\beta-\alpha|/H$ is close to $e$.  We use the curves
from 3.2:  $y^k = f(x)$ and the point $P = (\alpha, f(\alpha)^{1/k})$
and $[F(t)(P):F(t)] = kr$ now. Then $h(P) = ((d-1)(k-1)-2)H/kr$ by the
same calculation as before.  Also the ramification of $K(P)$ over
$K(\alpha)$ (where $K =F(t)$) is bounded by $(1+d-e)(k-1)H$.  So by
Riemann-Hurwitz $2g(P)-2 < k(2g(\alpha) -2) + (1+d-e)(k-1)H$, where
$g(\alpha)$ is the genus of $K(\alpha)/F$. To bound $g(\alpha)$ apply
the Castelnuovo inequality to $K(\alpha)$ viewed as the compositum of
$F(t)$ and $F(\alpha)$ both function fields of genus 0 together with
$[K(\alpha):F(t)] = r, [K(\alpha):F(\alpha)] = H$. So $g(\alpha) < Hr$
and $2g(P)-2 < 2krH + (1+d-e)(k-1)H$.  Finally, we can apply Theorem
2, provided the Kodaira-Spencer map is of maximal rank and get
$$((d-1)(k-1)-2)H/kr  < (2+\epsilon)(2krH + (1+d-e)(k-1)H)/kr + O(1).$$
Now divide by $H$ and make $H$ big and $\epsilon$ small, obtaining
$$((d-1)(k-1)-2)/kr  \le 2(2kr + (1+d-e)(k-1))/kr.$$
The last inequality gives a bound for $e$ in terms of $d,r$ and $k$.
If we can take $k$ arbitrarily large it gives
$e \le (d+3+4r)/2$.
Of course this is only interesting when $r < (d-3)/4$ but it seems that
other methods that yield improvements on Liouville's inequality, such as
Osgood's, do not give anything in this setting. 

4.0.2 In this case, our finite characteristic version of Vojta's 
conjecture gives (under the maximal rank condition) 
$e\leq 2r+2$, for $r>1$.  

\section{Explicit formulas}

5.0 Now we turn to calculation of Kodaira-Spencer matrix and explicit
formulae for the  rank conditions.

5.1 For $X$ as in 3.1, the basis for holomorphic differentials is
$\omega_{(a, b)}=x^ay^b (dx/F_y)$, with $a, b \geq 0$ and $a+b\leq
d-3$: Differentiating $F=1$ (treating $t$ as a constant for now), we
get $F_xdx+F_ydy=0$, hence $dx/F_y=-dy/F_x$. Since $F_x$ and $F_y$ are
not simultaneously zero, one has only to look at poles at
infinity. Since $dx$ has order two pole and $F_y$ has order $d-1$ pole
there, the claim follows.

Since $F_xdx+F_ydy+F_tdt=0$, the relative differential 
$dx/F_y$ has good liftings: $dx/F_y$ in the open set 
$U_1$ where $F_y\neq 0$ and $-dy/F_x$ in the open 
set $U_2$ where $F_x\neq 0$.  By our assumptions 
$U_1$ and $U_2$ cover $X$ and Cech cohomology computation 
shows that the map $KS: H^0(\Omega_{X/S}^1)\rightarrow
H^1({\cal O}_X)\otimes \Omega_S^1$ sends 
$dx/F_y$ to $F_tdt/(F_xF_y)$. Hence $x^ay^bdx/F_y$ is 
sent to $x^ay^bF_tdt/(F_xF_y)$.

We calculate the $g \times g$ matrix $M=(m_{ij})$ of 
the KS map in the basis above: 
If $P_k$'s are zeros of $F_y$, i.e points in the 
complement of $U_1$, then since the Serre duality 
sends a differential  on $U_1\cap U_2$ 
to the sum of its residues at $P_k$'s (if we use 
$U_2$ instead, we get negative of this: it is 
well-defined only up to a sign), we have 
$$m_{(a, b)(r, s)}=\sum_k Res_{P_k}(\frac{x^{a+r}y^{b+s}F_tdx}
{F_y^2F_x})dt \in \Omega_S^1$$

5.2 In the situation of 3.2, now $x$ has degree $k$ and $y$ has degree $d$, 
so that $x^idx/y^j$ is holomorphic, as long as $0<j<k$ and 
$(i+1)k+1\leq jd$. Note that $\sum_{j=1}^{k-1} \lfloor (jd-1)/k\rfloor
=\sum jd/k-\sum j/k = (k(k-1)/2)d/k-(k(k-1)/2)/k=g$, 
since $k$ and $d$ are relatively 
prime. Since  these differentials are linearly independent, they 
give the basis of the holomorphic differentials. 

Since $ky^{k-1}dy=f_xdx+f_tdt$, $dx/y^j$ has good liftings:
$dx/y^j$ in the open set $U_1$ where $y\neq 0$ and 
$ky^{k-1-j}dy/f_x$ in the open set $U_2$ where $f_x\neq 
0$. By our assumptions $U_1$ and $U_2$ cover $X$. The 
KS map sends $x^idx/y^j$ to $x^if_t dt/f_xy^j$. 

Similar calculation then gives
$$m_{(i, j), (l, n)}=\sum_s Res_{P_s}(\frac{x^{i+l}f_tdx}{y^{j+n}f_x})dt$$
where $P_s$ are now zeros of $y$, i.e. the conjugates 
of $\beta$. In the hyper-elliptic case $k=2$ (so that $d$ and $p$ are 
odd), this simplifies to 
$$m_{i,l}=\sum_s Res_{P_s}(\frac{x^{i+l}f_tdx}{y^2f_x})dt=
\sum_s \frac{P_s^{i+l}f_t(P_s)}{2f_x(P_s)^2}dt$$
since $y^2=f(x)$, so that $f/(x-P_s)|_{P_s}=f_x(P_s)$ 
and the fact that $x-P_s$ is of degree $2$ gives rise to 
the factor of $2$. 

\section{Remarks and questions}

6.0 If $\beta$ satisfies a rational Riccati equation,
then we know that the KS class of $Y:=$ the projective line minus the
Galois orbit of $\beta$ is zero and hence $Y$ is defined over
$F(t^p)$.  So for appropriate models, we have $f_t=F_t=0$ and since KS
is independent of co-ordinates and separable base change, we see that
KS is zero in this case for the examples in 3.1  and 3.2. In other
words, KS is non-zero implies `not Riccati' and hence our
inequalities, using Theorem 2, follow by Osgood's result proved under
weaker hypothesis.  Hence, the hierarchy given in 2.2, does not
give any new hierarchy in that case unconditionally
(except possibly in Wirsing-type result above
as well as in approximation results on non-rational base field that we
get using non-rational base $S$ in Theorem 2, where there are no
earlier results), though it suggests that there is such hierarchy
conjecturally, giving Vojta's inequality under the maximal rank. What
are the best inequalities one can conjecture?  (It is not just half
the bound, because that would be Roth for $g=2$, on just `not Riccati'
hypothesis, and that is known to be false by the examples in [V2]).
For the height inequalities, multiplier $2g-2$ would be best under
non-zero KS and $1$ would be best under the maximal rank. (So we
understand $g=2$ at least). What would be the best multipliers in between? 
Are they obtained by interpolating reciprocal-linearly as in 2.2? 

6.1 One can ask the similar question for the exponent hierarchy: 
But here different association of curves seem to lead to 
different conditions and bounds and the correct formulation is still 
unclear (even whether the hierarchy is finite or infinite),
except it is likely that the maximal rank (`generic') gives the 
exponent $2$ and KS non-zero would give some exponent between 
$d/2$ (attained by examples of [V2]) and Osgood bound. In 
this context, note that the condition that KS vanishes is independent 
of $k$. Is the maximal rank condition (or the whole hierarchy 
for that matter) also independent of $k$? Explicit calculation 
in 5.2 might help in deciding this.

6.2 It is well-known that general curves have maximal rank KS, but
though most of our curves (namely Thue curves for irreducible
polynomials and super-hyperellptic curves with branch points
consisting of Galois orbit together with infinity) are most probably
of maximal rank, since bounding a rank would give a differentially
closed condition, this has not been established. As a simple example,
if we are in characteristic $2$ and $\beta$ of degree $4$ is given by
$\beta^4+a\beta^3+b\beta^2+c\beta +d =0$, then KS is zero
(i.e. $\beta$ satisfies Riccati) implies $ac'=a'c$.  Now, if we can show
for any $d$, $k$ (sufficiently large), $p$ that there is at least one
curve $y^k=f(x)$, with $\deg f=d$ in characteristic $p$ with maximal
rank KS, then it is easy to show that most do: The coefficients of $f$
satisfy a differential equation, which can be turned into an algebraic
equation by writing each coefficient $a=\sum_0^{p-1}a_i^pt^i$, so
getting an algebraic equation on the $a_i$ which is not identically
zero so it is not satisfied by most $a_i$.  Our explicit calculations
of KS maps may help construct such examples, but it has not been done
yet.

6.3 We have optimistically suggested that Vojta inequality 
would hold in finite characteristic under the maximal KS 
rank hypothesis, but it may be that higher order deformation 
theory is needed for that. 

\section{Riccati and cross-ratios}

7.0 Finally, we record the proof of the claim  on 
Riccati connection (Note that vanishing of KS (i.e. 
having no infinitesimal deformations of first order)
is equivalent in this case to vanishing of derivatives 
of all cross ratios of 4 conjugates): Let $\beta\in F((1/t))$ be algebraic 
(and so automatically separable) of degree
$d$ over $F(t)$, with $\beta^d+b_{d-1}(t)\beta^{d-1}+\cdots 
=0$. Implicit differentiation gives 
$$\beta'=a_n\beta^n+\cdots +a_0,\  a_i\in F(t), \ 
a_n \neq 0, \ n\leq d-1 \ \ \ \ \ \ \ \ \ \ \ (1)$$

{\bf Claim}: $n=2$, i.e. $\beta$ satisfies the rational 
Riccati equation $\beta'=a\beta^2+b\beta+c$, with $a, b, c\in F(t)$ 
if and only if the cross ratio of any four conjugates 
of $\beta$ has zero derivative. 

{\bf Proof}: The derivative of the cross-ratio of 
$\beta, \beta_1, \beta_2, \beta_3$ being zero is equivalent to 
$$\frac{\beta'-\beta_1'}{\beta-\beta_1}+\frac{\beta_3'-\beta_2'}{\beta_3-\beta_2}
-\frac{\beta'-\beta_2'}{\beta-\beta_2}-\frac{\beta_3'-\beta_1'}{\beta_3-\beta_1}=0 \\ \ \ \ \ \ \ \ \  
\ \ \ \ \ \ (2)$$
Now conjugates $\beta_i$ of $\beta$ also satisfy (1). Hence

$$\frac {\beta_i'-\beta_j'}{\beta_i-\beta_j}=a_n\frac{\beta_i^n-\beta_j^n}{\beta_i-\beta_j}
+\cdots +a_2\frac{\beta_i^2-\beta_j^2}{\beta_i-\beta_j}+a_1$$

If $a_n=0$ for $n\geq 3$, the left hand side of (2) then reduces to
$a_2((\beta+\beta_1)+(\beta_3+\beta_2)-(\beta+\beta_2)-(\beta_3+\beta_1))$, which is zero. This proves
the `only if' statement. 

The `if' statement will be proved by
a contradiction: So we fix any 3 conjugates $\beta_1, \beta_2, \beta_3$ and 
assume that for any other conjugate $\beta$
we have 
$$\sum_{m>2}a_m (\frac{\beta^m-\beta_1^m}{\beta-\beta_1}+\frac{\beta_3^m-\beta_2^m}{\beta_3-\beta_2}
-\frac{\beta^m-\beta_2^m}{\beta-\beta_2}-\frac{\beta_3^m-\beta_1^m}{\beta_3-\beta_1})=0$$

Now $(\beta^m-\beta_1^m)/(\beta-\beta_1)=\sum_{i+j=m-1} \beta^i\beta_1^j$ and 
$\beta^j\beta_1^k+\beta_3^j\beta_2^k-\beta^j\beta_2^k-\beta_3^j\beta_1^k=(\beta-\beta_3)[\beta_1^k
(\beta^j-\beta_3^j)/(\beta-\beta_3)-\beta_2^k(\beta^j-\beta_3^j)/(\beta-\beta_3)]$, so that 
taking out the non-zero factor $(\beta-\beta_3)(\beta_1-\beta_2)$ we get 
$$0=\sum_{m>2}a_m[\sum_{j, k>0, j+k=m-1}
\frac{\beta^j-\beta_3^j}{\beta-\beta_3}\frac{\beta_1^k-\beta_2^k}{\beta_1-\beta_2}]$$
Now the quantity between the square brackets is just 
$\sum \beta_{i_1}\cdots \beta_{i_{m-3}}$, where each $\beta_i$ 
is one of the four conjugates. Hence we get 
$0=a_n\sum \beta_{i_1}\cdots \beta_{i_{n-3}}+a_{n-1}
\sum \beta_{i_1}\cdots \beta_{i_{n-4}}+\cdots + a_3$. 

The coefficient of $a_n$ is of degree $n-3$ in 
$\beta_i$'s. Subtracting a similar equation that 
one obtains when $\beta$ is replaced by another conjugate 
$\overline{\beta}$, and taking out the non-zero factor 
$\beta-\overline{\beta}$ we get another equation with degrees 
dropping by one. Continuing in this fashion with 
other conjugates (there are $d-4\geq n-3$ of them), 
we get $a_n=0$, a contradiction.

{\bf Acknowledgments}: The second author would like to thank 
Nitin Nitsure and Nicholas Katz for helpful conversations 
about the deformation theory. He would also like to thank 
Tata Institute of Fundamental Research where he visited 
while this research was carried out.

\end{document}